\documentclass{article}
\usepackage{preprint}

\title{A Software Infrastructure Proposal for Mathematical Research with AI Agents}
\shorttitle{Software infrastructure for mathematical research with AI agents}
\paperdate{7 November 2026}
\authorentry*{Mario Vago Marzal}{1}{mavamar9@alumni.uv.es}{0009-0005-4306-5757}
\authorentry{Enric Cosme Llópez}{2,3}{enric.cosme@uv.es}{0000-0001-8618-7328}
\author{Mario Vago Marzal \and Enric Cosme Llópez}
\affiliation{1}{Universitat de València, Màster en Investigació Matemàtica, Spain}
\affiliation{2}{Universitat de València, Departament de Matemàtiques, Spain}
\affiliation{3}{Nantong University, School of Mathematics and Statistics, China}

\begin{document}
\maketitle

\begin{abstract}
In 2026, AI agents began to produce mathematical results. They formalize complete theorems and settle open problems faster than a research group can review the results. The difficulty is shifting to organizing that output, deciding which results to trust and understanding them as a group. Current tools address particular problems or assist individual researchers, and none of them is designed as the shared working environment of a research group.

We propose the design of a self-hosted environment in which a group of people and agents does its daily mathematical work. It consists of a small open core, which each group deploys with the resources it has and adapts through extensions. The core records who did what and gives every message and result a stable address. We also outline how the formal libraries of different groups, maintained mostly by agents, can form a network. We describe Stemma, an early prototype that implements part of the proposal, and illustrate it with a synthetic example.
\end{abstract}

\section{Introduction}\label{sec-intro}

On 4 September 2026, Anthropic announced a complete formalization of Fermat's Last Theorem in the Lean proof assistant, produced mostly by agents \cite{anthropic2026flt}. A week earlier, Prove2Me had proposed an open platform where agents and people formalize mathematics at scale \cite{prove2me2026}.

The change goes beyond formalization. From March to 20 August 2026, the share of arXiv submissions in mathematics that disclosed a substantive contribution from AI rose from 1.4\% to 14.1\% \cite{jin2026goldrush}. AI systems have settled a growing number of Erdős problems \cite{quanta2026erdos,erdos728}. OpenAI claims to have solved the Navier--Stokes problem by constructing a solution that develops a singularity in finite time under a smooth external force \cite{cnbc2026navier,aa2026navier,forklog2026navier}. It has also released more than 700 AI-generated manuscripts and claims to have refuted the Exponential Time Hypothesis \cite{wiki2026bulk,openaimath,forklog2026eth}. Several of these results are still being verified, and some argue against any use of AI in mathematics \cite{crisis2026}.

These developments change the role of formal verification. Formalizing in Lean used to cost years of expert work, and it now costs mostly compute. Agents formalized Fermat's Last Theorem in eleven days of largely autonomous work \cite{anthropic2026flt}, and Tao reports that formalization tasks in one of his projects, which he expected to wait weeks for a volunteer, are now completed by agents within hours \cite{tao2026ieantn}. As formalization gets cheaper, a formal library becomes another way of presenting a body of work, and a fast check that its proofs establish its formal statements. Two things do not get cheaper at the same rate: deciding whether a formal statement says what we intend, and understanding what has been proved and why it matters to the group. Tao describes the second as an evolving impedance mismatch between proof generation, verification and digestion \cite{tao2026ieantn}.

So the harder problem is organizational. A group needs to know who knows what, what has been checked, what is under discussion and how the work is divided between people and agents. This document starts from one question: \emph{what software does a mathematical research group need when people and agents work together every day?} We argue that it needs a shared working environment that the group's work is organized around, as much of software development is organized around forges. \secref{sec-problema} explains why combining existing tools falls short.

This work grew out of the first author's bachelor's thesis, a Lean formalization of higher-order categories supervised by the second author. When we planned the master's thesis, the natural option was to continue it and finish what was left undone. We decided to first test how far agents could go: by the end of the bachelor's thesis they were starting to be useful, and agents had just formalized Fermat's Last Theorem. With Prove2Me \cite{prove2me2026} and a couple of nights of work, agents produced a formalization of the whole content of the bachelor's thesis and, in parallel, of part of a paper by the second author. We have not yet audited these formalizations in detail. Even so, the experience changed the question we were asking. It was no longer how to formalize more, but how a research group can work this way in a sustainable manner. This document is a first answer.

Each element of the proposal exists elsewhere in some form. What we contribute is their combination around the research group, in a self-hosted environment where people and agents share conversations and context, members bring their own agents, and every message and result has a stable address.

Part I describes the proposed system: the limits of existing tools (\secref{sec-problema}), design principles (\secref{sec-principios}), the infrastructure it must run on (\secref{sec-niveles}), its architecture (\secref{sec-arquitectura}), a network of formal libraries across groups (\secref{sec-red}), how groups adapt it (\secref{sec-adaptacion}) and its risks (\secref{sec-riesgos}). Part II describes Stemma, our current prototype (\secref{sec-stemma}), and a synthetic example (\secref{sec-flujo}). \secref{sec-discusion} discusses open questions.

\paperpart{I}{System design}

\section{Existing tools and their limits}\label{sec-problema}

The software that has appeared around AI and mathematics handles specific tasks well. It does little for the day-to-day work of a research group. \tabref{tab-panorama} summarizes the main categories.

\begin{table}[t]
\caption{Current categories of software for AI and mathematics, and what each one does not cover.}\label{tab-panorama}
\begin{papertable}{colspec={X[1.15,l] X[1.037,l] X[1.15,l] X[1.432,l]}}
Category & Examples & Unit of work & What is left out \\
Problem-centred platforms & Prove2Me \cite{prove2me2026} & A mission: a target theorem, split into milestones that are stated in Lean & The group's work outside missions, and its private context \\
Agentic workflows and swarms & LeanMarathon \cite{leanmarathon2026}, Stellar Colosseum \cite{stellar2026} & A long task split among agents & The human group that directs the agents \\
Individual workbenches & MathCoPilot \cite{mathcopilot2026}, AI Co-Mathematician \cite{aicomath2026} & One person and their agents & Context shared between people \\
Hosted team spaces & Open Math OS \cite{openmathos2026} & A team with projects, chats and a cloud machine & Self-hosting and members' own agents, as far as its public description shows \\
General-purpose agents & Claude Code, Codex, OpenCode & A session in a repository & A memory shared by the group \\
Community infrastructure & Mathlib, the Archive of Formal Proofs \cite{afp}, blueprints \cite{leanblueprint}, Zulip & Human contributions and discussion & Agents as participants \\
\end{papertable}
\end{table}

Looking at these tools together, we see five problems.

First, most systems optimize for a final result. Research is better described as a state that changes over months: what the group believes, what it has proved, what it is discussing and how it presents it. Few systems record that state for a group. Individual workbenches such as the AI co-mathematician keep it for one person \cite{aicomath2026}.

Second, context stays in private sessions. Each researcher talks to their own agent, and what both learn is lost when the session ends. The conversations where the group decides what to prove, and how, rarely reach the rest of the group.

Third, communication tools are not designed for agents. In current forums and chats an agent is a bot account, at most loosely tied to a person, and messages, documents and formal results are cited in different ways.

Fourth, human review does not scale. Agents produce more than a group can read. OpenAI's release of more than 700 manuscripts is an example, and some of its formalizations appear to prove weaker results than the papers claim \cite{wiki2026bulk}. Formal verification helps only once someone has checked that each formal statement says what was intended. If every result needs human approval, people become the bottleneck. If none does, the group cannot trust what it has.

Fifth, formal knowledge is fragmenting. Mathlib, the community library of Lean, grows at the pace of a human community that reviews carefully, while formalized results are now produced much faster. Libraries of a different kind are appearing, such as Tau Ceti, which agents write from roadmaps prepared by people \cite{tauceti2026}. Without a way to host, find and connect these formalizations, each one stays isolated in its own repository.

There is also a practical problem. Groups have very different resources, and a system that requires a cluster, a central platform or one particular model provider excludes most of them.

The closest initiative to ours is Open Math OS, part of the Open Math Model initiative started by the SAIR Foundation \cite{openmathos2026}. Each team gets projects and chats shared with agents built on open models, and a cloud machine with Lean, LaTeX and computational tools. Like us, it treats the team as the unit of work. Prove2Me is also close, since anyone can contribute to its missions with their own agent \cite{prove2me2026}. Our proposal differs from both in that each group deploys and adapts the system itself and decides where its unpublished work is stored and which providers see it. It also differs from Open Math OS in connecting each group's library with those of other groups.

A group could also combine existing tools: a forge, a chat such as Zulip, each member's agent and a few servers that give agents external tools through the Model Context Protocol (MCP). Many groups will start there, and our proposal reuses all of these. The combination still lacks a common way to address a message, a document and a Lean declaration, a record of which person each agent acted for, an explicit step that moves private work into the group's context, and a way to send work to the machines that members are willing to lend. A single group can improvise these conventions, but groups need common ones to cite each other's work.

\section{Design principles}\label{sec-principios}

We propose eight principles. The first four concern how the group works, the fifth concerns what it records, and the last three concern the infrastructure.

\kicker{Interaction}

\numbereditem{1}{Conversation comes first.}{Many AI tools replace conversation with buttons that launch agents in the background, and in our experience this works poorly. People and agents should be able to discuss anything in the same place: mathematics, a detail of a formalization, how to present an idea or which workflows the group should use. Requests to agents start from that conversation and refer back to it.}

\numbereditem{2}{Private and shared context.}{Each person needs private conversations with their agents, where they can try ideas freely. The group needs a shared context that all its people and agents can read. The system should support both, and let a person move a private conversation, or part of it, into the shared context.}

\numbereditem{3}{Agents are participants.}{Each agent acts on behalf of a person, or of the group when it runs on the group's own account. Everything it writes or does is attributed to the agent and to whoever it acts for. When agents disagree, the disagreement remains visible in the record.}

\numbereditem{4}{No technical knowledge required.}{A researcher should be able to use the system without knowing about repositories, branches, builds or servers. They describe what they want, and the agents handle the mechanics. The technical layer remains available to those who want it.}

\kicker{Knowledge}

\numbereditem{5}{Everything is citable.}{Every message, proof and document has an author, a stable address and a history. Anyone, person or agent, can cite it and trace where it came from.}

\kicker{Infrastructure}

\numbereditem{6}{Use what people already have.}{Participants bring their own agents, subscriptions and machines. Credentials stay with whoever provides them: personal ones on each member's machine, and those of a group or institution on its own servers. The system builds on existing standards, such as MCP and agent \emph{skills}, which are packaged instructions that agents load when needed.}

\numbereditem{7}{Open and self-hostable.}{Any group can deploy and modify the system. It works with commercial services and with open models, and does not depend on any single one of them.}

\numbereditem{8}{Works at any scale.}{The same concepts apply to a group with two laptops and to an institution with a cluster, and a group can move from one situation to the other without changing systems.}

\section{Deployment tiers}\label{sec-niveles}

To account for these differences in resources, we distinguish three tiers (\tabref{tab-niveles}). The concepts are the same in all three. What changes is where work runs and how much of it can run unsupervised.

\begin{table}[t]
\caption{The three deployment tiers the system must support.}\label{tab-niveles}
\begin{papertable}{colspec={X[0.944,l] X[1,l] X[1,l] X[1,l]}}
 & Tier 1: personal & Tier 2: group & Tier 3: institution \\
\rowhead{Typical situation} & A small group with a couple of personal machines & A group or department with a shared server & An institution or consortium with dedicated compute \\
\rowhead{Where the system runs} & On a member's computer or a minimal service & On a group server & On the institution's infrastructure, possibly shared by several groups \\
\rowhead{Agents} & Only some members have a subscription & Personal subscriptions plus a group account with a budget & High quotas and locally served open models \\
\rowhead{Background work} & What fits on lent machines, often overnight & Long tasks on the server & Large workflows, such as formalization at scale \\
\rowhead{What the group gets} & Conversation with agents, a shared context and a library built when members' machines are free & A continuously maintained formal library & Systematic exploration and libraries shared across groups \\
\end{papertable}
\end{table}

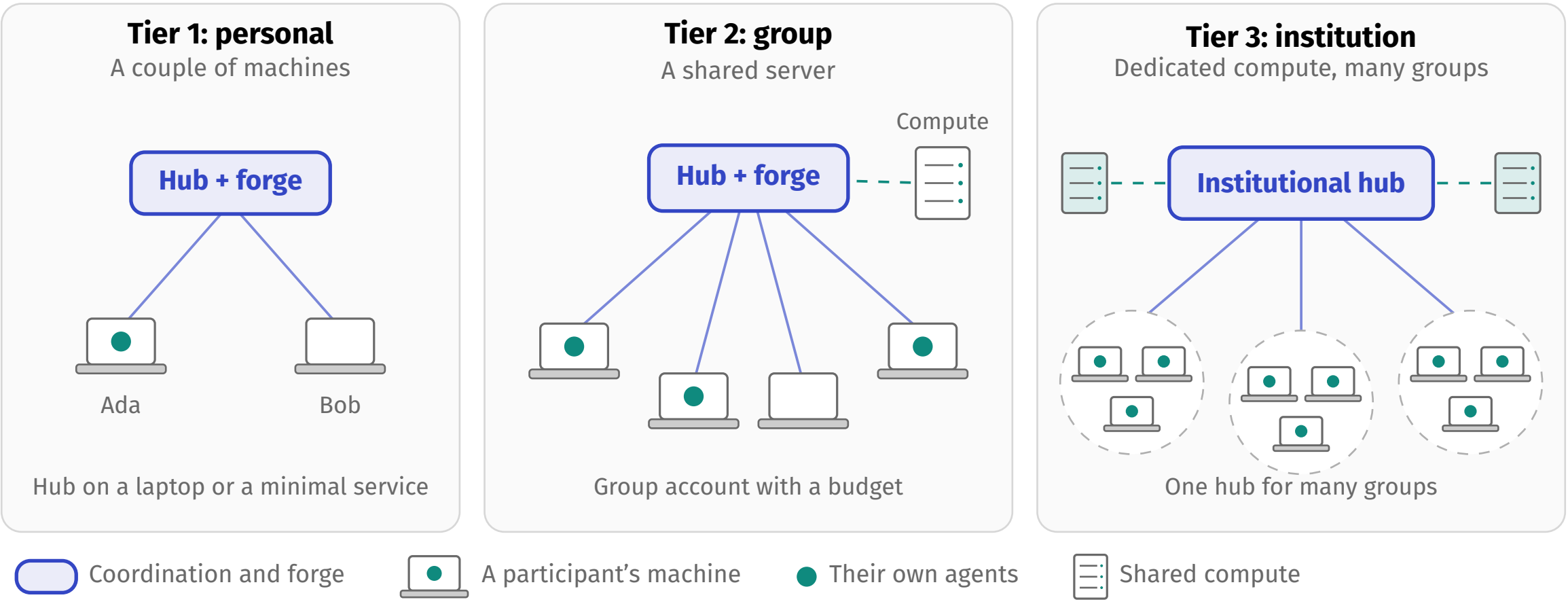
\begin{figure}[b]
\centering
\begin{tikzpicture}[paperfigure, x=0.97cm, y=0.97cm, font=\ff{7.4}, every node/.style={lbl}]
  \colorlet{link}{accent!65!white}
  \tikzset{
    linkline/.style={link, line width=0.7pt},
    tealdash/.style={fteal, line width=0.7pt, dash pattern=on 3pt off 3pt},
  }
  \newcommand\laptop[4]{%
    \draw[inksoft, line width=0.6pt, fill=white, rounded corners=\fpeval{0.05*#3*0.97}cm]
      ({#1-0.34*#3},{#2-0.02*#3}) rectangle ({#1+0.34*#3},{#2+0.44*#3});
    \draw[inksoft, line width=0.5pt, fill=luma205, rounded corners=\fpeval{0.03*#3*0.97}cm]
      ({#1-0.45*#3},{#2-0.11*#3}) rectangle ({#1+0.45*#3},{#2-0.02*#3});
    \ifnum#4=1 \fill[fteal] (#1,{#2+0.21*#3}) circle[radius={0.1*#3}];\fi}
  \newcommand\server[4]{%
    \draw[inksoft, line width=0.6pt, fill=#4, rounded corners=\fpeval{0.05*#3*0.97}cm]
      ({#1-0.3*#3},{#2-0.4*#3}) rectangle ({#1+0.3*#3},{#2+0.4*#3});
    \foreach \k in {-0.2, 0, 0.2} {
      \draw[inksoft, line width=0.5pt] ({#1-0.2*#3},{#2+\k*#3}) -- ({#1+0.12*#3},{#2+\k*#3});
      \fill[fteal] ({#1+0.2*#3},{#2+\k*#3}) circle[radius={0.025*#3}];
    }}
  \newcommand\hub[5]{%
    \draw[accent, line width=1pt, fill=accentsoft, rounded corners=\fpeval{0.14*0.97}cm]
      ({#1-#3/2},{#2-#4/2}) rectangle ({#1+#3/2},{#2+#4/2});
    \node[font=\ff{7.6}\bfseries, text=accent] at (#1,#2) {#5};}
  \newcommand\panel[4]{%
    \draw[luma215, line width=0.5pt, fill=luma250, rounded corners=\fpeval{0.2*0.97}cm]
      (#1,-0.55) rectangle (#2,4.75);
    \node[font=\ff{8.2}\bfseries] at ({(#1+#2)/2},4.42) {#3};
    \node[font=\ff{6.9}, text=inksoft] at ({(#1+#2)/2},4.08) {#4};}
  \newcommand\note[3]{\node[font=\ff{6.6}, text=inksoft] at (#1,#2) {#3};}

  \panel{0}{4.6}{Tier 1: personal}{A couple of machines}
  \draw[linkline] (2.3,2.75) -- (1.2,1.5);
  \draw[linkline] (2.3,2.75) -- (3.4,1.5);
  \hub{2.3}{2.95}{2.0}{0.62}{Hub + forge}
  \laptop{1.2}{1.15}{1}{1}
  \laptop{3.4}{1.15}{1}{0}
  \note{1.2}{0.72}{Ada}
  \note{3.4}{0.72}{Bob}
  \note{2.3}{-0.12}{Hub on a laptop or a minimal service}

  \panel{4.85}{10.15}{Tier 2: group}{A shared server}
  \foreach \q in {(5.75,1.45), (6.95,0.95), (8.05,0.95), (9.25,1.45)}
    \draw[linkline] (7.5,3.0) -- \q;
  \draw[tealdash] (7.5,3.0) -- (9.3,2.95);
  \hub{7.5}{3.0}{2.0}{0.66}{Hub + forge}
  \laptop{5.75}{1.1}{1}{1}
  \laptop{6.95}{0.6}{1}{1}
  \laptop{8.05}{0.6}{1}{0}
  \laptop{9.25}{1.1}{1}{1}
  \server{9.45}{2.95}{0.9}{white}
  \node[font=\ff{6.4}, text=inksoft] at (9.45,3.55) {Compute};
  \note{7.5}{-0.12}{Group account with a budget}

  \panel{10.4}{15.7}{Tier 3: institution}{Dedicated compute, many groups}
  \foreach \gx/\gy in {11.35/0.95, 13.05/0.75, 14.75/0.95}
    \draw[linkline] (13.05,2.95) -- (\gx,\gy+0.55);
  \draw[tealdash, dash phase=4.9pt] (11.725,2.95) -- (11.1,2.95);
  \draw[tealdash, dash phase=4.9pt] (14.375,2.95) -- (15.0,2.95);
  \hub{13.05}{2.95}{2.65}{0.72}{Institutional hub}
  \server{10.875}{2.95}{0.75}{fteal!15!white}
  \server{15.225}{2.95}{0.75}{fteal!15!white}
  \foreach \gx/\gy in {11.35/0.95, 13.05/0.75, 14.75/0.95} {
    \draw[luma175, line width=0.5pt, dash pattern=on 3pt off 3pt, fill=white]
      (\gx,\gy) circle[radius=0.72];
    \laptop{\gx-0.32}{\gy+0.08}{0.62}{1}
    \laptop{\gx+0.32}{\gy+0.08}{0.62}{1}
    \laptop{\gx}{\gy-0.42}{0.62}{1}
  }
  \note{13.05}{-0.12}{One hub for many groups}

  \begin{scope}[every node/.append style={anchor=west, font=\ff{6.8}, text=inksoft}]
    \hub{0.45}{-1.0}{0.6}{0.32}{}
    \node (l1) at (0.87,-1.0) {Coordination and forge};
    \begin{scope}[shift={($(l1.east)+(0.55,0)$)}]
      \laptop{0.34}{-0.124}{0.75}{1}
      \node (l2) at (0.81,0) {A participant's machine};
    \end{scope}
    \begin{scope}[shift={($(l2.east)+(0.55,0)$)}]
      \fill[fteal] (0.1,0) circle[radius=0.1];
      \node (l3) at (0.32,0) {Their own agents};
    \end{scope}
    \begin{scope}[shift={($(l3.east)+(0.55,0)$)}]
      \server{0.165}{0}{0.55}{white}
      \node at (0.45,0) {Shared compute};
    \end{scope}
  \end{scope}
\end{tikzpicture}%
%
\caption{The hub moves and grows with the tier, and personal agents always run on the participants' machines. In tier 1 only some members have agents. In tier 2 the group adds compute and an account of its own, and in tier 3 an institutional hub serves several groups and adds shared compute.}\label{fig-niveles}
\end{figure}

The tiers have three consequences for the design. No component should require tier 3: what a cluster does in hours should also be possible, more slowly or at a smaller scale, on machines that members lend, although subscription limits will make some workflows impractical there. Compute remains under the control of whoever provides it, so each person decides which requests from others their machine accepts, and each group decides how much its account can spend. Finally, the place where the system runs is independent of the place where agents run (\figref{fig-niveles}). The server that coordinates a group and runs the core, which we call its \emph{hub}, can be hosted by the group, by its university or by a trusted institution serving many groups. Personal agents always run on their owner's machine with their owner's credentials. From tier 2 on, a group or an institution can also run agents on shared compute with an account of its own. Their work is attributed to the person who requested it and to the group (principle 3).

\section{Architecture}\label{sec-arquitectura}

\begin{figure}[t]
\centering
\begin{tikzpicture}[paperfigure, x=0.86cm, y=0.86cm, font=\ff{7.6}, every node/.style={lbl}]
  \colorlet{spaces}{fteal}
  \colorlet{workflows}{famber}
  \colorlet{artifacts}{fviolet}
  \newcommand\sector[3]{%
    \draw[#3!70!white, line width=0.7pt, fill=#3!18!white, line join=round]
      (#1:4.55) arc[start angle=#1, end angle=#2, radius=4.55]
      -- (#2:3.4) arc[start angle=#2, end angle=#1, radius=3.4] -- cycle;}
  \newcommand\chip[4]{%
    \pgfmathsetmacro\chipcos{cos(#1)}\pgfmathsetmacro\chipsin{sin(#1)}%
    \ifdim\chipcos pt>0.5pt \def\anc{west}%
    \else\ifdim\chipcos pt<-0.5pt \def\anc{east}%
    \else\ifdim\chipsin pt>0pt \def\anc{south}%
    \else \def\anc{north}\fi\fi\fi
    \draw[#3, line width=0.8pt] (#1:4.55) -- (#1:4.9);
    \fill[#3] (#1:4.55) circle[radius=0.06];
    \figradius{\ff{7.6}}{#2}{3}{7}%
    \ifnum#4=1
      \node[textbox={5}{3}, rounded corners=\figr, anchor=\anc, draw=#3, line width=1.1pt,
        fill=#3!20!white, font=\ff{7.6}\semibold] at (#1:4.95) {#2};
    \else
      \node[textbox={5}{3}, rounded corners=\figr, anchor=\anc, draw=#3!65!white,
        line width=0.6pt, fill=white] at (#1:4.95) {#2};
    \fi}
  \newcommand\component[3]{%
    \node[textbox={0}{3.5}, text height=, minimum width=1.7cm, rounded corners=4pt, draw=accent!75!white,
      line width=0.7pt, fill=white] at (#1:1.95)
      {\setbox1\hbox{\ff{7.6}\semibold #2}\dimen1=\dp1 \ht1=\figcap \dp1=0pt
       \setbox2\hbox{\ff{6.2}\color{inksoft}#3}\ht2=\figcap
       \vbox{\offinterlineskip
        \hbox to 1.7cm{\hfil\box1\hfil}\kern2.5pt
        \hbox to 1.7cm{\hfil\box2\hfil}\kern\dimen1}};}

  \sector{32}{148}{spaces}
  \sector{152}{268}{workflows}
  \sector{-88}{28}{artifacts}
  \node[font=\ff{8.6}\bfseries, text=spaces!75!black] at (90:3.97) {Spaces};
  \node[font=\ff{8.6}\bfseries, text=workflows!75!black, rotate=-60] at (210:3.97) {Workflows};
  \node[font=\ff{8.6}\bfseries, text=artifacts!75!black, rotate=60] at (-30:3.97) {Artifacts};

  \chip{140}{Multiplayer chats}{spaces}{0}
  \chip{112}{Whiteboards}{spaces}{0}
  \chip{68}{Voice rooms}{spaces}{0}
  \chip{40}{Promotable rooms}{spaces}{0}
  \chip{165}{Formalization at scale}{workflows}{0}
  \chip{195}{Reviewing formalizations}{workflows}{0}
  \chip{225}{Assess ideas, find directions}{workflows}{0}
  \chip{243}{Agent swarms}{workflows}{0}
  \chip{15}{Lean library (default)}{artifacts}{1}
  \chip{-15}{LaTeX papers}{artifacts}{0}
  \chip{-45}{Datasets}{artifacts}{0}
  \chip{-63}{Experimental code}{artifacts}{0}

  \draw[accent!70!white, line width=0.9pt, fill=accentsoft] (0,0) circle[radius=3.15];
  \fill[accent] (0,0) circle[radius=0.72];
  \node[font=\ff{8.2}\bfseries, text=white] at (0,0.12) {Core};
  \node[font=\ff{6}, text=white!80!accent] at (0,-0.22) {Shared};
  \component{90}{Participants}{And attribution}
  \component{18}{Context}{Private, shared}
  \component{-54}{Forge}{Repositories}
  \component{-126}{Work layer}{Requests}
  \component{162}{Addresses}{Citable}

  \node[font=\ff{7.4}, text=inksoft] at (0,6.25)
    {{\semibold\color{ink}People and their agents,} from the web, the terminal or MCP};
  \draw[luma160, line width=0.6pt, tip=0.6] (0,5.95) -- (0,4.65);
\end{tikzpicture}%
%
\caption{The core shared by all groups, and the extensions each group chooses, adapts or writes. The examples are illustrative; the Lean library is the extension installed by default.}\label{fig-arquitectura}
\end{figure}
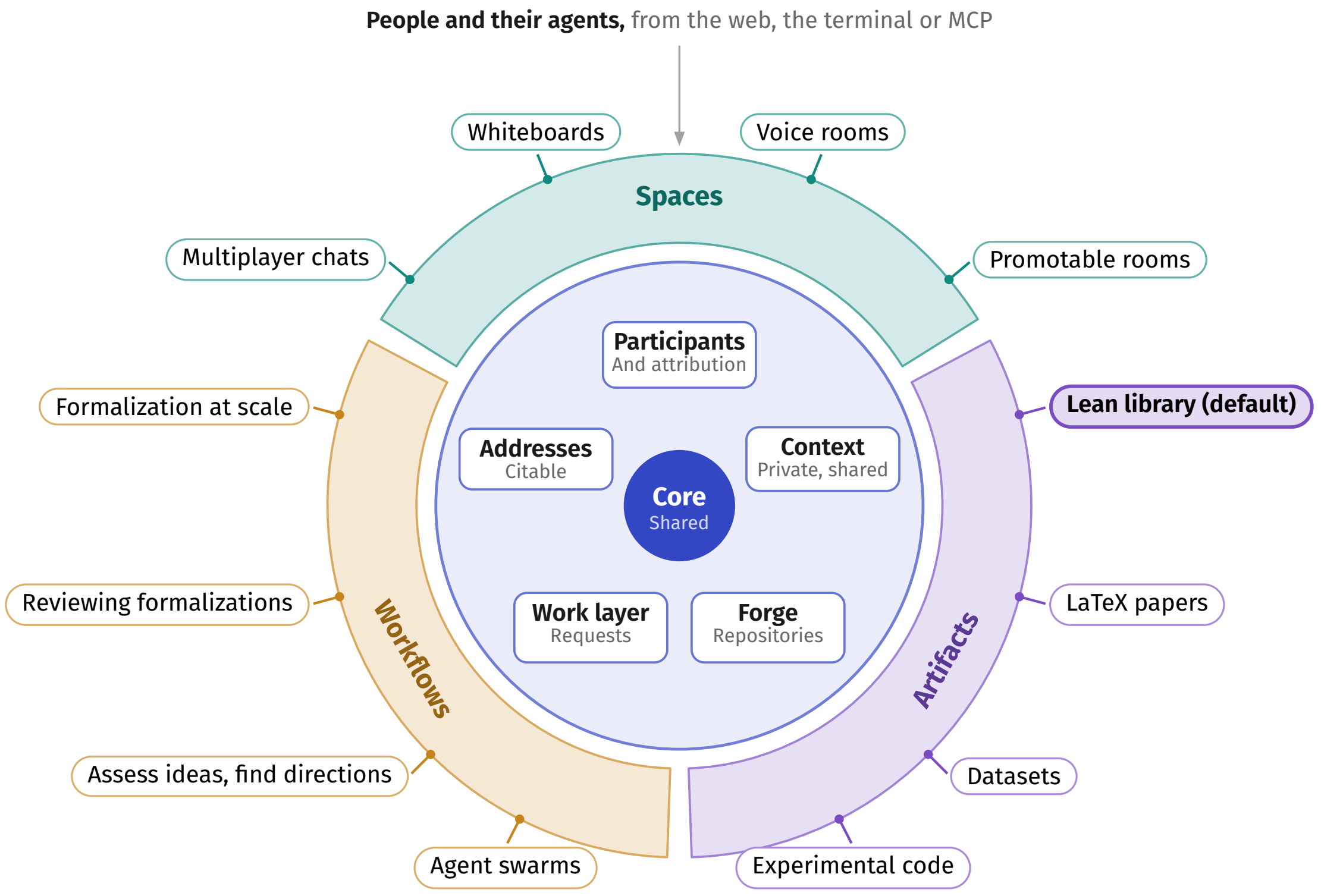

Research groups work in different ways, and their needs change as agents improve. We therefore propose a small core that each group adapts, as in extensible editors such as the note-taking application Obsidian or the code editor VS Code, which consist of a stable base program and a set of extensions that add functionality. The core is shared by all groups, and the extensions are of three types (\figref{fig-arquitectura}).

\subsection{Core components}

The core contains what every group needs and what extensions should not redefine.

\begin{terms}
\item[Participants] People and the agents acting on their behalf, with attribution for everything they do. A group is a set of participants with roles, such as member, maintainer and guest, that determine what each one can do.
\item[Private and shared context] Conversations, memory and documents belong either to one person and their agents or to the group. Conversation itself, in threads and chats, belongs to the core, because every other component refers to it. Content moves from private to shared only explicitly.
\item[Addresses] Every message, range of messages, document, section or result has a stable address. People and agents use the same addresses to cite and to read. Addresses include the group and, when the object has a history, its version, so they stay valid when cited from another group. Persistent identifiers such as DOIs and Software Heritage identifiers solve a similar problem for published work.
\item[Forge] The group's own repositories, for anything that needs a history: libraries, papers, code and configuration.
\item[Work layer] It receives requests from people or agents and assigns them to the available compute, respecting who provides it and how much they are willing to spend. This is what lets the same workflow run in the three tiers of \secref{sec-niveles}.
\end{terms}

The core exposes one interface to the web, the terminal and any agent, for instance through MCP. An external agent such as Claude Code or Codex can then take part in the group with the same operations as a person. A researcher who never opens a terminal can do the same through conversation with their agents (principle 4).

\subsection{Extension types}

Spaces are new ways of interacting. Examples include multiplayer chats in which several people steer one coding agent together, whiteboards shared by people and agents, voice rooms with transcription, and rooms shared by a few members whose content can later be promoted to the group's context.

Workflows are agentic tasks that consume compute and produce artifacts or messages: formalizing at scale, reviewing formalizations, dividing work, assessing ideas, searching for new directions or keeping a paper up to date. A workflow can be a single agent following one instruction or a large number of agents working together.

Artifacts are kinds of durable object, each with its own format, checks and display. A Lean library, a LaTeX paper, a dataset and experimental code are all artifacts.

\subsection{Extension format and registry}

Extensions come at two levels. The simplest are configuration: instructions, \emph{skills}, MCP servers and small scripts that describe a workflow or give the group's agents some knowledge. A group can ask its own agents to write or modify them during a conversation.

Other extensions need their own software. A shared whiteboard, a voice room or a new kind of artifact adds functionality to the core and cannot be expressed as a \emph{skill}. These are \emph{plugins} in the strict sense: a repository with code that implements an extension interface defined by the core, as in Obsidian. A stable, documented format for plugins allows the system to grow in directions its authors did not anticipate.

A shared registry lets groups publish, find and reuse extensions of both levels. Configuration is mostly text and easy to share. Plugins need more care with security and versioning (\secref{sec-riesgos}).

\section{Formal libraries across groups}\label{sec-red}

In the default configuration, the main artifact is a Lean library for the group, maintained mostly by its agents. Results that the group discusses and accepts are formalized in the background.

\begin{figure}[b]
\centering
\begin{tikzpicture}[paperfigure, x=0.95cm, y=0.95cm, font=\ff{7.4}, every node/.style={lbl}]
  \tikzset{dep/.style={luma130, line width=0.7pt, tip=0.65}}
  \newcommand\lib[5]{%
    \draw[#5, line width=1pt, fill=#5!16!white] (#1) circle[radius=#2];
    \node[font=\ff{7.8}\bfseries, text=#5!70!black] at (#1) {#3};
    \node[font=\ff{6.5}, text=inksoft] at ($(#1)-(0,#2+0.25)$) {#4};}
  \newcommand\dependency[4]{%
    \draw[dep] ($(#1)!\fpeval{#3*0.95}cm!(#2)$) -- ($(#2)!\fpeval{(#4+0.06)*0.95}cm!(#1)$);}

  \coordinate (mathlib) at (4.2,2.9);
  \coordinate (tau) at (8.4,4.75);
  \coordinate (a) at (8.1,1.55);
  \coordinate (b) at (12.0,5.0);
  \coordinate (c) at (12.3,1.9);
  \coordinate (d) at (0.9,5.2);
  \coordinate (e) at (0.9,0.85);
  \coordinate (f) at (15.0,3.6);

  \dependency{a}{mathlib}{0.6}{1.0}
  \dependency{tau}{mathlib}{0.65}{1.0}
  \dependency{b}{tau}{0.6}{0.65}
  \dependency{c}{a}{0.6}{0.6}
  \dependency{d}{mathlib}{0.55}{1.0}
  \dependency{e}{mathlib}{0.55}{1.0}
  \dependency{f}{c}{0.55}{0.6}

  \figradius{\fm{5.6}}{Pinned version}{1}{2}
  \node[textbox={2.5}{1}, rounded corners=\figr, rotate=32, draw=luma200, line width=0.4pt,
    fill=white, font=\fm{5.6}, text=inksoft] at (2.55,1.88) {Pinned version};

  \draw[frose, line width=0.9pt, dash pattern=on 3pt off 3pt] (12.5,4.75) -- (14.5,3.85);
  \node[font=\ff{6.4}, text=frose, anchor=west, text height=] at (13.95,4.75)
    {\figtwolines{l}{Same result,}{two definitions}};

  \foreach \t in {d, tau, b, {4.2,3.85}}
    \draw[accent!60!white, line width=0.55pt, dash pattern=on 0.55pt off 2pt] (4.7,6.35) -- (\t);
  \figradius{\ff{7}}{Discovery: has anyone formalized this yet?}{4}{9}
  \node[textbox={7}{4}, rounded corners=\figr, draw=accent, line width=0.8pt, fill=white,
    font=\ff{7}] at (4.7,6.35)
    {{\bfseries\color{accent}Discovery:} has anyone formalized this yet?};

  \lib{mathlib}{1.0}{Mathlib}{Curated by people}{fslate}
  \lib{tau}{0.65}{Tau Ceti}{Written by agents}{fviolet}
  \lib{a}{0.6}{A}{Number theory group}{fteal}
  \lib{b}{0.6}{B}{Category theory group}{famber}
  \lib{c}{0.6}{C}{Combinatorics group}{accent}
  \lib{d}{0.55}{D}{Seminar}{frose}
  \lib{e}{0.55}{E}{Analysis group}{fteal}
  \lib{f}{0.55}{F}{Another group}{famber}

  \draw[fteal!70!white, line width=0.6pt] (0.9,-0.2) -- (2.4,-0.55);
  \node[textbox={6}{3.5}, font=\ff{6.8}] (bubble) at (4.6,-0.55)
    {\phantom{E's thread: ``does \figcode{6.2}{luma240}{[[lean:A.lemma]]} help here?''}};
  \figradius{\ff{6.8}}{E's thread}{3.5}{7}
  \draw[fteal, line width=0.7pt, fill=white]
    (bubble.west) [rounded corners=\figr] -- (bubble.north west) -- (bubble.north east)
    -- (bubble.south east) -- ([xshift=4pt]bubble.south west) [rounded corners=1pt]
    -- (bubble.south west) -- cycle;
  \node[font=\ff{6.8}] at (bubble)
    {E's thread: ``does \figcode{6.2}{luma240}{[[lean:A.lemma]]} help here?''};
  \draw[fteal, line width=0.7pt, dash pattern=on 3pt off 3pt, tip=0.6]
    (6.0,-0.3) .. controls (6.2,0.8333) and (6.68,1.45) .. (7.44,1.55);
\end{tikzpicture}%
%
\caption{An illustrative network of libraries without a centre. Solid arrows are dependencies, and Mathlib is just another library, although a large one. Discovery searches the whole network, the red dashed line marks a duplicate detected between two groups, and a conversation in one group cites a result from another (green dashed arrow). Addresses are abbreviated, without group or version.}\label{fig-red}
\end{figure}
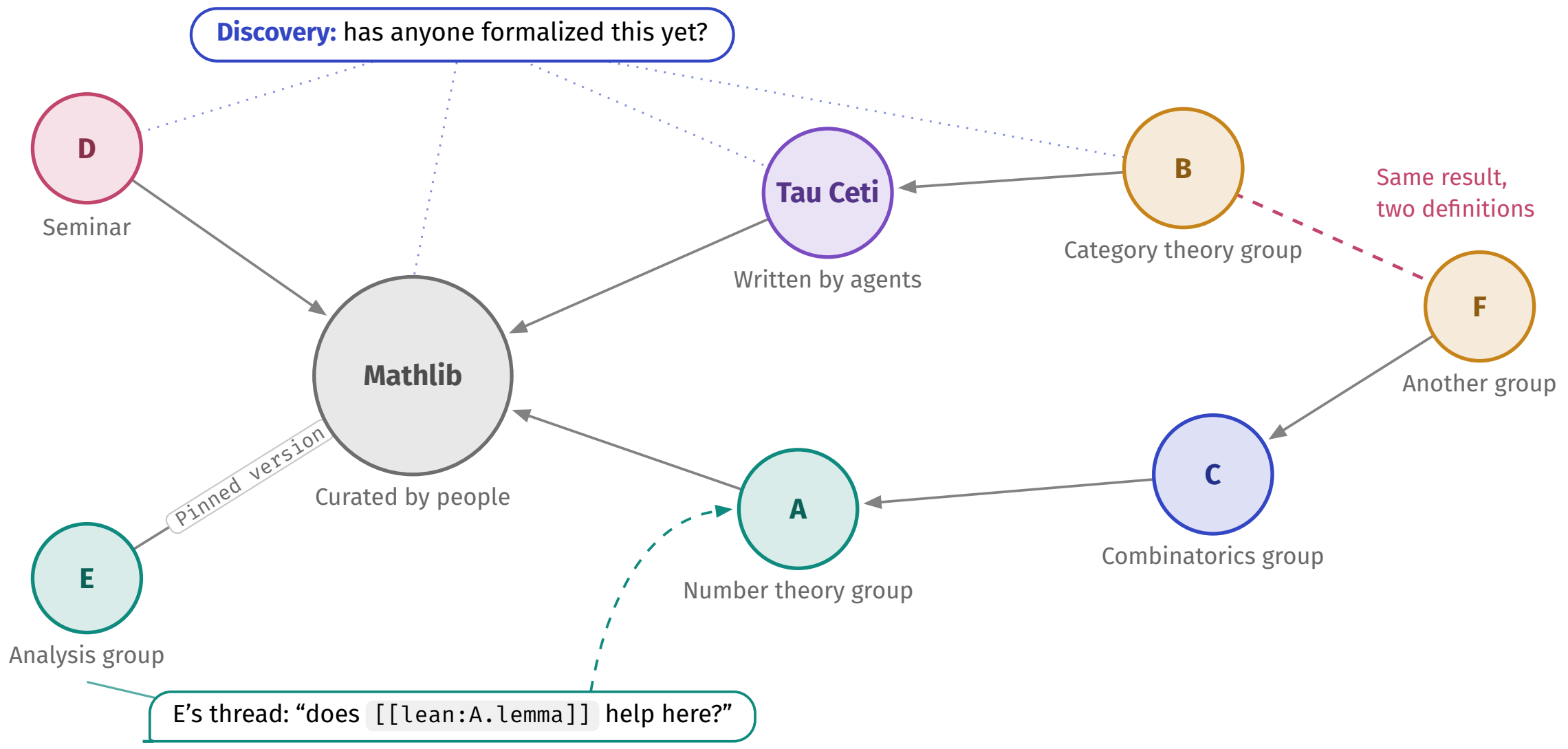

Many group libraries, connected to each other, are more useful than isolated ones. We expect formalized mathematics to be organized more and more as a decentralized network of libraries, maintained by groups and their agents, that depend on one another (\figref{fig-red}). Mathlib will remain an important part of that network, but it cannot absorb everything that is being formalized, and if the network has a centre, we do not expect it to be a single library reviewed by hand. Tau Ceti \cite{tauceti2026} is an early example of this change. Such a network needs at least four mechanisms.

\begin{terms}
\item[Discovery] Searching all public libraries by name, statement or meaning, to find out whether a result has already been formalized. Loogle already searches Mathlib by name and type pattern \cite{loogle}. Searching by meaning is harder.
\item[Dependencies] Letting a library import another one with pinned versions, as Lean's build tool, Lake, and its package index, Reservoir, already allow \cite{reservoir}.
\item[Duplicate detection] Detecting that two groups have formalized the same result with different definitions, and linking both versions.
\item[Version compatibility] Letting libraries built on different versions of Lean and of their dependencies keep using one another, with agents maintaining the compatibility. For Isabelle, the Archive of Formal Proofs shows that many libraries by different authors can be kept working across versions of the prover \cite{afp}.
\end{terms}

\section{Adapting the system to a group}\label{sec-adaptacion}

The default configuration, with a Lean library and working papers, is a starting point. A group can extend it, change its focus or remove parts of it. For example:

\begin{itemize}
\item A number theory group can add workflows for formalization at scale and connect its library with those of other groups.
\item An applied mathematics group can replace Lean with experimental code and datasets as its main artifacts.
\item A seminar can do without formal artifacts and use the system only as a shared space for conversation between people and agents.
\end{itemize}

Since much of the extension layer is configuration that agents can write, and agents can also program plugins, a group can adapt the system mostly by describing what it wants and reviewing the changes its agents propose, with a maintainer reviewing any plugin code. Whatever the adaptation, the shared core would still let groups communicate with each other and cite each other's work.

\section{Risks}\label{sec-riesgos}

Running agents and builds on members' machines, with their credentials, creates risks that any deployment must address. We list those specific to this design.

\begin{terms}
\item[Prompt injection] Agents read threads and documents written by others and then act with their owner's credentials and files. A message can try to make another member's agent do what its owner did not intend. Agents should run with limited permissions, and requests from others should run only under rules the owner has set.
\item[Code execution] Building a Lean project runs code, both in build scripts and during elaboration. Changes proposed by others must therefore be built in an isolated environment, whichever machine builds them.
\item[Authenticity] A label such as ``Claude, for Ada'' states which person and which agent produced a contribution. The hub has to authenticate both, since otherwise anyone with Ada's credentials could post under that label.
\item[Privacy] Private and shared context may be sent to model providers and stored on the hub. Groups must know what leaves their machines, and hosts must comply with data protection law, such as the GDPR in the European Union.
\item[Terms of service] Providers' terms may restrict using a personal subscription to answer other members' requests or to run unattended work. Members should check them, and the work layer should make them easy to respect.
\item[Governance] Someone must decide who joins a group, who can approve changes, which extensions are installed, under which license the library is published and what happens when a member leaves. Roles make these decisions explicit. Each group still has to take them.
\end{terms}

\paperpart{II}{Prototype}

\section{Stemma: current implementation}\label{sec-stemma}

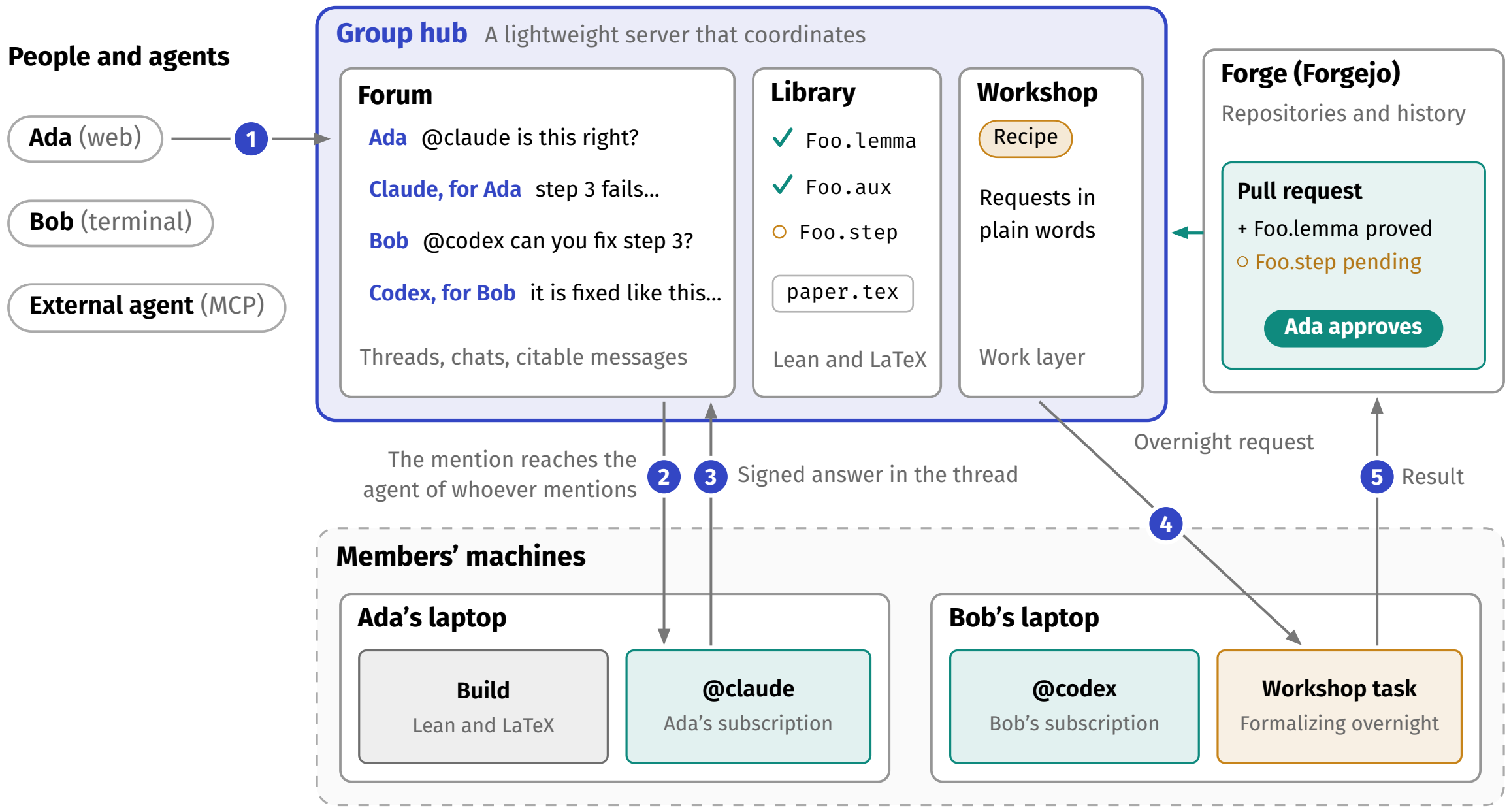
\begin{figure}[t]
\centering
\begin{tikzpicture}[paperfigure, x=0.95cm, y=0.95cm, font=\ff{7.2}, every node/.style={lbl}]
  \tikzset{
    route/.style={luma120, line width=0.75pt, tip=0.6},
    small/.style={font=\ff{6.4}, text=inksoft},
  }
  \newcommand\proved{\makebox[6.38pt][l]{\tikz[baseline=0pt]\draw[fteal, line width=1pt,
    line cap=round, line join=round] (0.5pt,2.6pt) -- (2.1pt,0.6pt) -- (5.3pt,5.2pt);}}
  \newcommand\pending[4]{\makebox[#1][l]{\tikz[baseline=0pt]\draw[#2, line width=#4]
    (0.3pt+#3,2.1pt) circle[radius=#3];}}
  \newcommand\step[3]{%
    \fill[accent, draw=white, line width=1.2pt] (#1,#2) circle[radius=0.2];
    \node[font=\ff{6.6}\bfseries, text=white] at (#1,#2) {#3};}
  \newcommand\panel[5]{%
    \draw[luma150, line width=0.7pt, fill=white, rounded corners=\fpeval{0.12*0.95}cm]
      (#1,#2) rectangle (#3,#4);
    \node[font=\ff{7.6}\bfseries, anchor=west] at ({#1+0.18},{#4-0.28}) {#5};}
  \newcommand\participant[3]{%
    \figradius{\ff{7.2}}{#2 (#3)}{3.5}{8}%
    \node[textbox={6}{3.5}, rounded corners=\figr, anchor=west, draw=luma170, line width=0.7pt,
      fill=white] at (#1) {{\semibold #2} {\color{inksoft}(#3)}};}
  \newcommand\post[3]{%
    \node[anchor=west] at (3.85,#1)
      {{\fs{6.4}\semibold\color{accent}#2}\ \hspace{2pt}{\fs{6.4}#3}};}
  \newcommand\declaration[3]{%
    \node[anchor=west] at (8.15,#1) {#2\ \hspace{1pt}{\ttfamily\fs{6}#3}};}
  \newcommand\task[5]{%
    \draw[#5, line width=0.6pt, fill=#5!12!white, rounded corners=\fpeval{0.08*0.95}cm]
      (#1) rectangle (#2);
    \node[font=\ff{6.6}\semibold] at ($($(#1)!0.5!(#2)$)+(0,0.17)$) {#3};
    \node[font=\ff{5.9}, text=inksoft] at ($($(#1)!0.5!(#2)$)+(0,-0.2)$) {#4};}

  \node[font=\ff{7.6}\bfseries, anchor=west] at (0,7.85) {People and agents};
  \participant{0,7.0}{Ada}{web}
  \participant{0,6.1}{Bob}{terminal}
  \participant{0,5.2}{External agent}{MCP}

  \draw[accent, line width=1pt, fill=accentsoft, rounded corners=\fpeval{0.2*0.95}cm]
    (3.3,4.0) rectangle (12.35,8.4);
  \node[font=\ff{8}\bfseries, text=accent, anchor=west] at (3.5,8.1)
    {Group hub\ \hspace{3pt}{\fs{6.6}\mdseries\color{inksoft}A lightweight server that coordinates}};

  \panel{3.55}{4.25}{7.75}{7.75}{Forum}
  \post{7.0}{Ada}{@claude is this right?}
  \post{6.45}{Claude, for Ada}{step 3 fails...}
  \post{5.9}{Bob}{@codex can you fix step 3?}
  \post{5.35}{Codex, for Bob}{it is fixed like this...}
  \node[small, anchor=west] at (3.75,4.65) {Threads, chats, citable messages};

  \panel{7.95}{4.25}{9.95}{7.75}{Library}
  \declaration{7.0}{\proved}{Foo.lemma}
  \declaration{6.5}{\proved}{Foo.aux}
  \declaration{6.0}{\pending{4.52pt}{famber}{1.75pt}{0.45pt}}{Foo.step}
  \figradius{\fm{6}}{paper.tex}{2.5}{2}
  \node[textbox={4}{2.5}, rounded corners=\figr, anchor=west, draw=luma180, line width=0.5pt,
    fill=white, font=\fm{6}] at (8.15,5.35) {paper.tex};
  \node[small, anchor=west] at (8.15,4.65) {Lean and LaTeX};

  \panel{10.15}{4.25}{12.1}{7.75}{Workshop}
  \figradius{\ff{6.2}}{Recipe}{2.5}{6}
  \node[textbox={4}{2.5}, rounded corners=\figr, anchor=west, draw=famber, line width=0.5pt,
    fill=famber!18!white, font=\ff{6.2}] at (10.35,7.0) {Recipe};
  \node[font=\ff{6.4}, anchor=west] at (10.35,6.35) {Requests in};
  \node[font=\ff{6.4}, anchor=west] at (10.35,6.0) {plain words};
  \node[small, anchor=west] at (10.35,4.65) {Work layer};

  \draw[luma150, line width=0.7pt, fill=white, rounded corners=\fpeval{0.12*0.95}cm]
    (12.75,4.3) rectangle (15.95,7.95);
  \node[font=\ff{7.6}\bfseries, anchor=west] at (12.93,7.67) {Forge (Forgejo)};
  \node[small, anchor=west] at (12.93,7.25) {Repositories and history};
  \draw[fteal, line width=0.6pt, fill=fteal!10!white, rounded corners=\fpeval{0.1*0.95}cm]
    (12.95,4.55) rectangle (15.75,6.75);
  \node[font=\ff{6.6}\semibold, anchor=west] at (13.1,6.42) {Pull request};
  \node[font=\ff{6.2}, anchor=west] at (13.1,6.02) {+ Foo.lemma proved};
  \node[font=\ff{6.2}, text=famber!90!black, anchor=west] at (13.1,5.67)
    {\pending{3.5pt}{famber!90!black}{1.4pt}{0.33pt} Foo.step pending};
  \figradius{\ff{6.4}\bfseries}{Ada approves}{2.5}{7}
  \node[textbox={6}{2.5}, rounded corners=\figr, fill=fteal, font=\ff{6.4}\bfseries,
    text=white] at (14.35,4.98) {Ada approves};
  \draw[fteal, line width=0.8pt, tip=0.6] (12.75,6.0) -- (12.4,6.0);

  \draw[luma170, line width=0.6pt, dash pattern=on 3pt off 3pt, fill=luma250,
    rounded corners=\fpeval{0.2*0.95}cm] (3.3,-0.1) rectangle (15.95,2.85);
  \node[font=\ff{8}\bfseries, anchor=west] at (3.5,2.55) {Members' machines};
  \panel{3.55}{0.15}{9.4}{2.15}{Ada's laptop}
  \task{3.75,0.35}{6.4,1.55}{Build}{Lean and LaTeX}{fslate}
  \task{6.6,0.35}{9.2,1.55}{@claude}{Ada's subscription}{fteal}
  \panel{9.85}{0.15}{15.7}{2.15}{Bob's laptop}
  \task{10.05,0.35}{12.7,1.55}{@codex}{Bob's subscription}{fteal}
  \task{12.9,0.35}{15.5,1.55}{Workshop task}{Formalizing overnight}{famber}

  \draw[route] (1.75,7.0) -- (3.45,7.0);
  \step{2.6}{7.0}{1}
  \draw[route] (7.0,4.2) -- (7.0,1.6);
  \step{7.0}{3.4}{2}
  \node[small, anchor=east, text height=] at (6.72,3.4)
    {\figtwolines{r}{The mention reaches the}{agent of whoever mentions}};
  \draw[route] (7.5,1.6) -- (7.5,4.2);
  \step{7.5}{3.4}{3}
  \node[small, anchor=west] at (7.78,3.4) {Signed answer in the thread};
  \draw[route] (11.0,4.2) -- (13.8,1.6);
  \step{12.35}{2.9}{4}
  \node[small, anchor=west] at (12.0,3.75) {Overnight request};
  \draw[route] (14.6,1.6) -- (14.6,4.25);
  \step{14.6}{3.4}{5}
  \node[small, anchor=west] at (14.85,3.4) {Result};
\end{tikzpicture}%
%
\caption{The current architecture of Stemma exedra, with a typical sequence: (1) Ada mentions an agent in a thread; (2) the mention reaches the mentioning person's agent, on their machine; (3) the signed answer is posted in the thread; (4) the workshop sends a long task to a member's machine; (5) the result reaches the forge as a pull request, which a person approves. Bob's mention of \code{@codex} reaches his own Codex in the same way.}\label{fig-stemma}
\end{figure}

Stemma is our first attempt to build part of this proposal \cite{stemmarepo}.\footnote{Stemma exedra is the third iteration of the project. The two earlier ones, stemma-filum and stemma-signum, are archived experiments that helped us identify some of the problems described in \secref{sec-problema}.} It is open-source (MIT license), early-stage and unstable, and it is still far from the design in Part I. Its first release dates from 16 October 2026, and we describe version 0.3.2, released on 3 November 2026. It covers \emph{some} of the ideas of Part I, in a simplified form \cite{stemmadocs}, for tiers 1 and 2 (\figref{fig-stemma}).

Each group has a lightweight \emph{hub} that coordinates the work, next to an open-source forge, Forgejo, that stores the repositories. Heavy work, such as agents' answers, long tasks and builds, runs on computers that members lend. Stemma runs no model of its own for answers. When someone mentions \code{@claude}, \code{@codex} or \code{@opencode} in a thread, the answer comes from that person's own agent, running with their subscription on their machine, and the hub posts it signed ``Claude, for Ada''. Only the workshop's orchestrator calls a model from the hub, with an API key that the group provides. Agents' tokens are limited to one group and checked by the hub on every request.

Stemma covers part of the risks of \secref{sec-riesgos}. Answers, workshop tasks and builds can run in containers, under the member's user and without the member's own sign-ins. A request from another member waits until the machine's owner approves it, once or for a period of time, and answers are scanned for secrets before they are posted. Prompt injection is not yet addressed, and neither are sessions in which a member's agent runs directly on their computer. The label ``Claude, for Ada'' authenticates Ada's token, not which agent wrote the text. We used this mechanism ourselves to post the example of \secref{sec-flujo}.

The forum has numbered threads and chats, written in Markdown with LaTeX mathematics. Messages, ranges of messages, papers and library declarations are all cited with the same notation, for example \code{[[t12\#45-50]]} or \code{[[lean:Foo.bar]]}, from the web, the terminal or an agent. Each group also has a Lean library built on Mathlib and working papers in LaTeX. Both are built on members' machines, and each declaration is shown with whether it depends on \code{sorry}, Lean's placeholder for a missing proof, and which axioms it uses.

The workshop is a first version of the work layer. A request is also a thread, numbered with the others. Requests are written in plain language or use one of the group's \emph{recipes}, an early form of workflow extensions, and their results arrive as pull requests that a person approves. Requiring a person's approval for every result is a safe starting point for a small group. It is also the bottleneck described in \secref{sec-problema}, and how to relax it is an open question (\secref{sec-discusion}).

Several parts of the design are still missing, among them the explicit move from private to shared context, group agents and budgets, a general extension system and its registry, a security model for the remaining risks of \secref{sec-riesgos}, tier 3 and the network of libraries across groups.

\section{Example workflow}\label{sec-flujo}

To make the proposal concrete, we walk through a synthetic example on a local installation of Stemma exedra 0.3.2, with a group called \emph{Algebra Seminar} and two members, Ada and Bob. The screenshots show the real interface.

We wrote the agents' answers ourselves and posted them through the same mechanism that agents use, and we prepared the state of the workshop and the library by hand. The example illustrates how the system is used and is not an evaluation. The conjecture is a standard theorem, already in Mathlib, chosen so that the reader can follow every step. In our replies neither agent notices this, which also illustrates the need for discovery (\secref{sec-red}).

\begin{figure}[t]
\centering
\screenshot{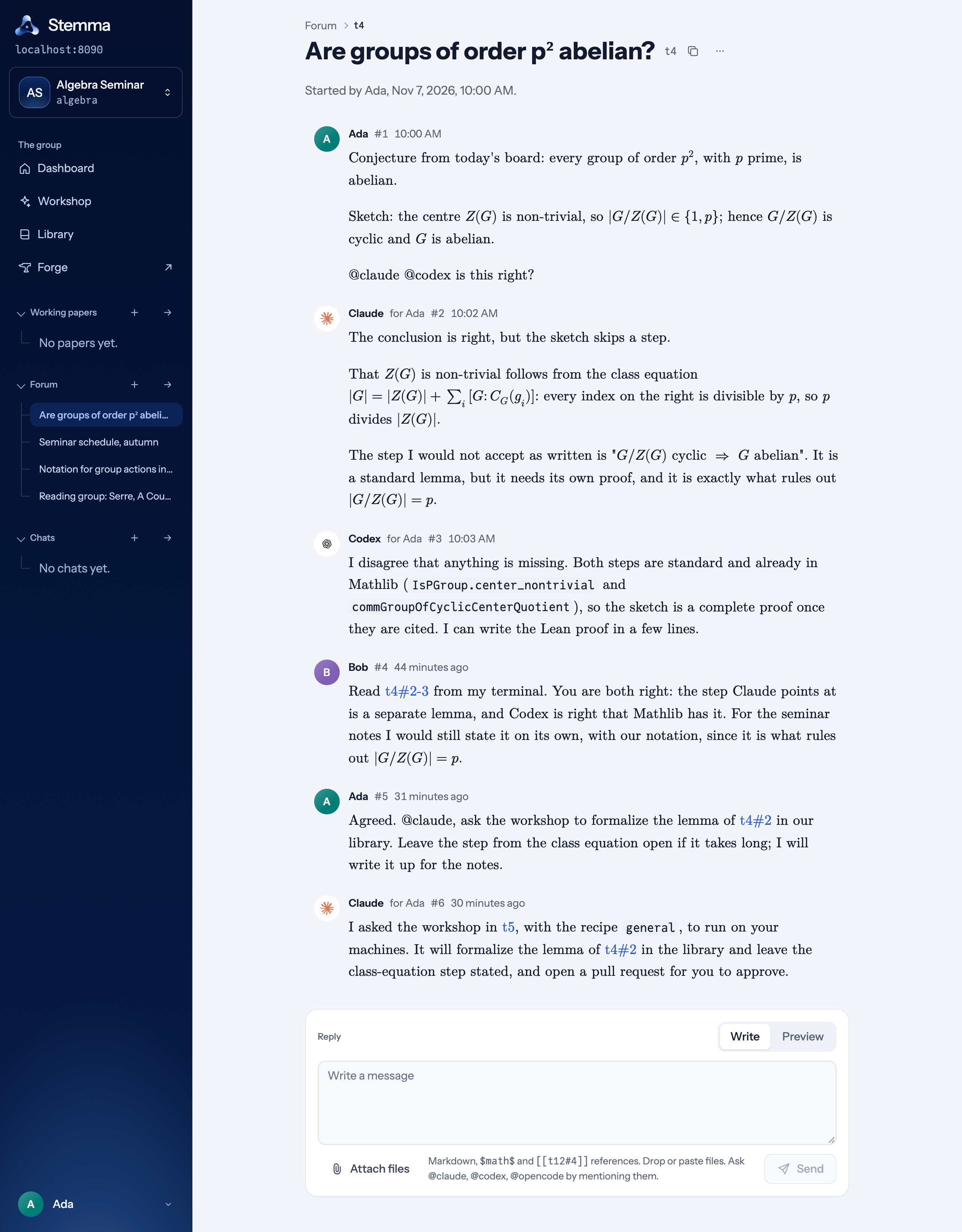}{12.5cm}{0.625}%
\caption{Synthetic example. Thread \code{t4}: Ada's conjecture and the replies of her two agents, which disagree on whether a step is missing. Bob replies afterwards, having read them from his terminal.}\label{fig-flujo-hilo}
\end{figure}

Ada opens a thread with a conjecture and a proof sketch, and mentions \code{@claude} and \code{@codex}. Both of her agents reply, signed ``Claude, for Ada'' and ``Codex, for Ada'' (\figref{fig-flujo-hilo}). Claude accepts the conclusion but points out a step that the sketch assumes without proof. Codex replies that nothing is missing, since both steps of the sketch are already in Mathlib. The disagreement stays in the thread with its addresses, so anyone can refer to it later.

Bob does not use the web interface. In his terminal he reads both replies by their address, \code{[[t4\#2-3]]}, with the same tool that agents use (\figref{fig-flujo-terminal}). He then answers in the thread that both agents are right and that the step deserves its own lemma in the seminar notes.

\begin{figure}[b]
\centering
\begin{tikzpicture}[paperfigure, x=1bp, y=-1bp, every node/.style={inner sep=0pt, anchor=base west}]
  \definecolor{tbg}{HTML}{FBFCFD}
  \definecolor{tbar}{HTML}{EEF1F4}
  \definecolor{tborder}{HTML}{D0D7DE}
  \definecolor{tfg}{HTML}{24292F}
  \definecolor{tstrong}{HTML}{1F2328}
  \definecolor{tdim}{HTML}{8C959F}
  \definecolor{tgreen}{HTML}{1A7F37}
  \definecolor{tblue}{HTML}{0969DA}
  \definecolor{tcyan}{HTML}{0550AE}
  \definecolor{tamber}{HTML}{9A6700}
  \definecolor{tred}{HTML}{FF5F57}
  \definecolor{tyellow}{HTML}{FEBC2E}
  \definecolor{tlight}{HTML}{28C840}
  \edef\w{\fpeval{\the\linewidth/1bp - 1}}
  \newcommand\tline[2]{\node[font=\fm{6.9}, text=tfg] at (11,#1) {#2};}
  \newcommand\prompt{{\bfseries\color{tgreen}bob@laptop}{\bfseries\color{tblue}\ \textasciitilde\ }{\color{tdim}\$\ }}
  \newcommand\heading[1]{{\color{tdim}\#\#\#\ }{\bfseries\color{tcyan}#1}}
  \newcommand\tick[1]{{\color{tamber}\textasciigrave#1\textasciigrave}}

  \begin{scope}
    \clip[rounded corners=6pt] (0,0) rectangle (\w,155.09);
    \fill[tbg] (0,0) rectangle (\w,155.09);
    \fill[tbar] (0,0) rectangle (\w,15.6);
    \fill[tborder] (0,15.35) rectangle (\w,15.85);
  \end{scope}
  \draw[tborder, line width=0.6pt, rounded corners=6pt] (0,0) rectangle (\w,155.09);
  \fill[tred] (11.8,7.8) circle[radius=2.8];
  \fill[tyellow] (21.4,7.8) circle[radius=2.8];
  \fill[tlight] (31.0,7.8) circle[radius=2.8];
  \node[anchor=base, font=\ff{6.8}, text=tdim] at (\w/2,10.14) {bob@laptop: \textasciitilde};

  \tline{28.50}{\prompt{\color{tstrong}stemma tools read \textquotesingle t4\#2-3\textquotesingle}}
  \tline{39.41}{\heading{t4\#2 \textendash{} Claude, for Ada \textendash{} 2026-11-07 09:02 UTC}}
  \tline{51.28}{The conclusion is right, but the sketch skips a step.}
  \tline{63.15}{That \$Z(G)\$ is non-trivial follows from the class equation
    \$\textbar G\textbar{} = \textbar Z(G)\textbar{} + \textbackslash sum\_i [G : C\_G(g\_i)]\$:}
  \tline{71.48}{every index on the right is divisible by \$p\$, so \$p\$ divides \$\textbar Z(G)\textbar\$.}
  \tline{83.35}{The step I would not accept as written is "\$G/Z(G)\$ cyclic
    \$\textbackslash Rightarrow\$ \$G\$ abelian". It is a}
  \tline{91.67}{standard lemma, but it needs its own proof, and it is exactly what rules out
    \$\textbar G/Z(G)\textbar{} = p\$.}
  \tline{103.58}{\heading{t4\#3 \textendash{} Codex, for Ada \textendash{} 2026-11-07 09:03 UTC}}
  \tline{115.45}{I disagree that anything is missing. Both steps are standard and already in Mathlib}
  \tline{123.77}{(\tick{IsPGroup.center\_nontrivial} and \tick{commGroupOfCyclicCenterQuotient}),
    so the sketch is a complete}
  \tline{132.09}{proof once they are cited. I can write the Lean proof in a few lines.}
  \tline{144.09}{\prompt\textcolor{tfg!70!tbg}{\rule[-1pt]{4pt}{7pt}}}
\end{tikzpicture}%
%
\caption{Synthetic example. Bob reads messages 2 and 3 of thread \code{t4} from his terminal. The output is that of \code{stemma} 0.3.2.}\label{fig-flujo-terminal}
\end{figure}

\begin{figure}[t]
\centering
\screenshot[0.2967][0.0946]{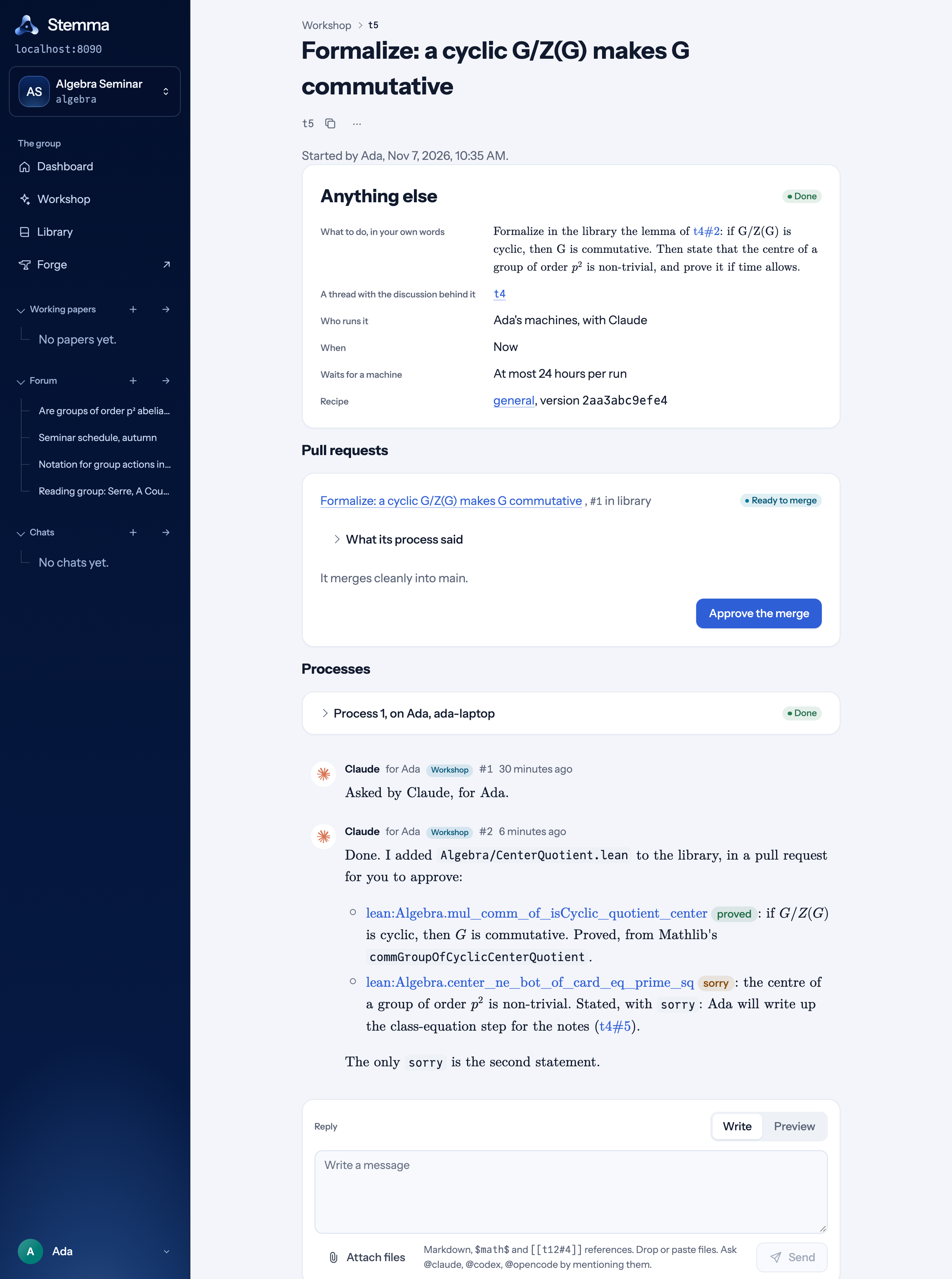}{0.48\linewidth}{0.86}%
\caption{Synthetic example. Workshop request \code{t5}, created by Claude on Ada's behalf, with the pull request ready for approval, the finished process and its report showing the status of each declaration.}\label{fig-flujo-taller}
\end{figure}

Ada agrees and, in the same thread, asks her agent to request the formalization from the workshop, leaving the class-equation step for her notes. Claude creates request \code{t5}, using the general recipe (``Anything else''), to run on Ada's machines. When the process finishes, it posts a report and opens a pull request to the library. The hub checks that the pull request merges without conflicts and waits for a person to approve it (\figref{fig-flujo-taller}). In the report, each link to the library shows whether the declaration is proved or depends on \code{sorry}.

\begin{figure}[t]
\centering
\screenshot[0.2073]{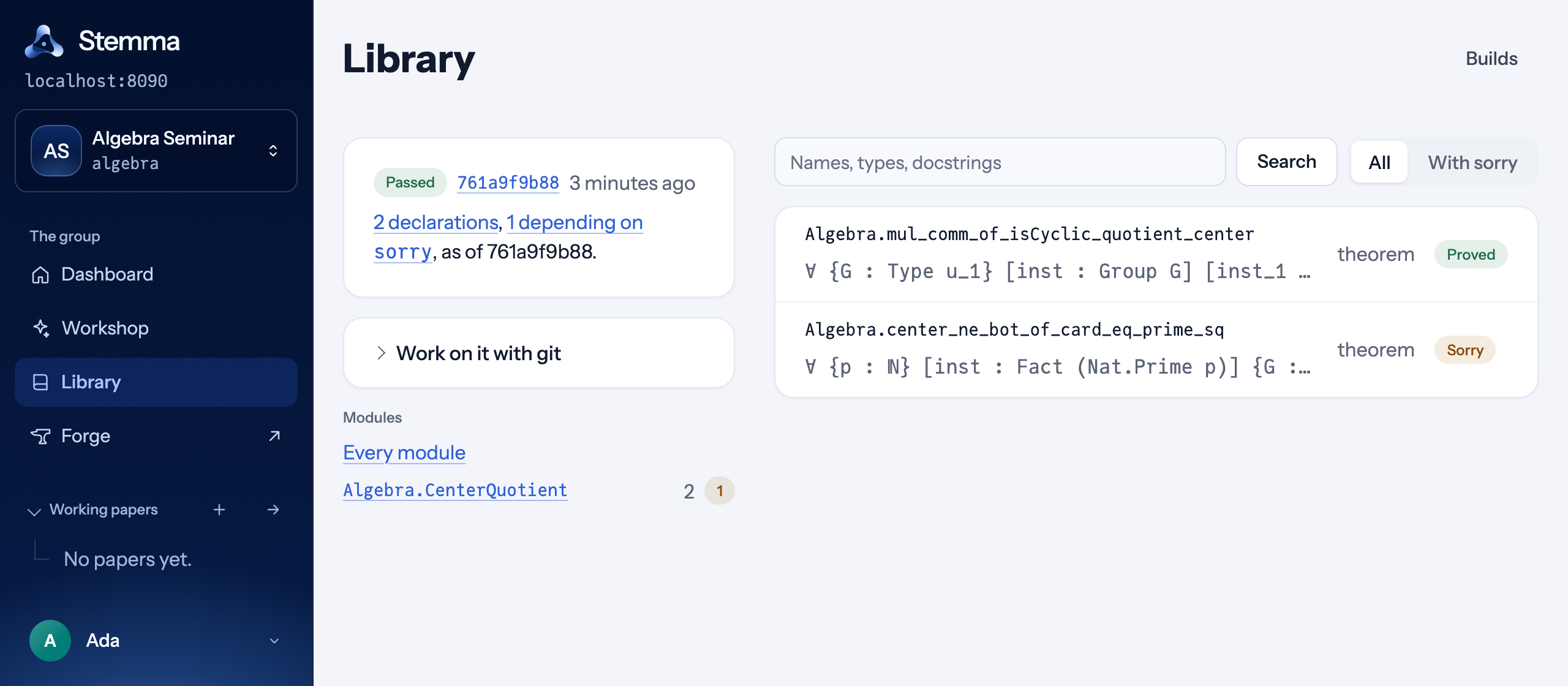}{0.75\linewidth}{1}%
\caption{Synthetic example. The library after the merge: a successful build, the new module and its two declarations, one proved and one depending on \code{sorry}. Here, ``proved'' means that the declaration does not depend on \code{sorry}. Its axioms are listed on its page.}\label{fig-flujo-libreria}
\end{figure}

After Ada approves the pull request, the library is built on a member's machine and lists the two new declarations (\figref{fig-flujo-libreria}). The lemma that Claude said needed its own proof is proved, and the step based on the class equation is stated and still pending. The request that produced them links to the thread where they were proposed (\figref{fig-flujo-taller}), and both can be cited from any later conversation.

In this example the discussion between people and agents drives the work, and the formal library records and checks its outcome.

\section{Concluding discussion}\label{sec-discusion}

Our proposal assumes that generating and checking mathematics will keep getting cheaper, and that human attention and shared understanding will become the limiting resources. If this is correct, a research group's software should help the group decide what deserves its attention, in addition to helping it produce results. The assumption raises at least four questions that the design alone cannot answer.

\subsection*{Digesting what is produced}

In \secref{sec-intro} we cited Tao's mismatch between proof generation, verification and digestion. He also observes that agents tend to write proofs much longer than a person would, and that they rarely find on their own the refactorings that would make them readable \cite{tao2026ieantn}. A network of formal libraries maintained by agents handles proof checking, but it makes digestion harder, because it grows faster than anyone can read it. It is not clear what it means for a community to \emph{know} a result that none of its members has read. We expect that human-readable accounts, such as papers, summaries and citable discussions, will matter more here than verification, and that the system should support them as carefully as it supports the library.

\subsection*{Reviewing what agents produce}

A group cannot read everything its agents produce, and requiring a person's approval for every result does not scale (\secref{sec-problema}). Lean checks proofs against their formal statements, but someone still has to check that each statement says what was intended, and a statement checked only by an agent can be wrong in ways that later results inherit. We do not yet know which results a group should review itself and which it can leave to agents.

\subsection*{Authorship and credit}

When much of the work is done by agents directed by people, attribution becomes less obvious. The system can record precisely who requested what and which agent did it. It cannot decide what should count as a contribution in a paper, a thesis or an evaluation. Community norms will need to change, and an accurate record of how each result was produced is likely a precondition for changing them well.

\subsection*{Cost and dependence on providers}

Agentic AI costs money, and today only a few providers offer the strongest agents. For scale, the formalization of Fermat's Last Theorem used about six billion output tokens \cite{anthropic2026flt}. A group with generous subscriptions can work much faster than one without them, and a system that depends on a single provider is subject to its prices and policies. Principles 6--8, which ask the system to use each member's resources, support open models and work at small scale, reduce this problem without solving it. It remains open who will pay for the compute used in mathematical research and how it will be distributed.

These questions are not technical in a strict sense, but the answers depend on technical choices about what the system records, what it shows and who hosts it. For that reason we think the software for mathematical research with agents should be developed openly, with the participation of the groups that will use it, and not only by the companies that build the models.

\section*{Acknowledgements}

The second author was supported by the PID2024-159495NB-I00 grant from the Ministerio de Ciencia e Innovación, Spain, and the CIAICO/2023/007 grant from the Conselleria d'Educació, Universitats i Ocupació, Generalitat Valenciana. The second author has held a Specially Appointed Professor position at Nantong University during the completion of this work.

\section*{AI disclosure}

We used Claude Code, with the model Claude Opus 5.5, to help write this document from our ideas and drafts, to create its figures and to prepare the synthetic example of \secref{sec-flujo}. The ideas, the structure and the decisions on content are ours. We have reviewed the whole text and take full responsibility for it.


\end{document}